\documentclass[a4paper,12pt,final]{article}
\usepackage{amsmath, amssymb,latexsym, amsthm}
\usepackage{graphicx}
\usepackage{epsfig}

\usepackage{amsmath,amsfonts,amssymb,amsthm}
\newtheorem{thm}{Theorem}
\newtheorem*{thm*}{Theorem}

\newtheorem*{example*}{Example}
\newtheorem{Def}{Definition}
\newtheorem*{Def*}{Definition}

\newcommand{\reftit}{\textit}    
\newcommand{\refis}{\textbf}     

\title{Adiabatic times for Markov chains and applications}
\author{Kyle Bradford\footnote{Department of Mathematics,  Oregon State University, Corvallis, OR  97331-4605, USA \texttt{bradfork@math.oregonstate.edu}}  
\quad and \quad
Yevgeniy Kovchegov \footnote{Department of Mathematics,  Oregon State University, Corvallis, OR  97331-4605, USA \texttt{kovchegy@math.oregonstate.edu}}}
\date{}

\begin{document}
\maketitle

\abstract{We state and prove a generalized adiabatic theorem for Markov chains and provide examples and applications related to  Glauber dynamics of Ising model over $\mathbb{Z}^d/n\mathbb{Z}^d$. The theorems derived in this paper describe a type of adiabatic dynamics for $\ell^1(\mathbb{R}^n_+)$ norm preserving, time inhomogeneous Markov transformations, while quantum adiabatic theorems deal with $\ell^2(\mathbb{C}^n)$ norm preserving ones, i.e. gradually changing unitary dynamics in $\mathbb{C}^n$. }


\section*{Introduction}
The long-term stability of time inhomogeneous Markov processes is an active area of research in the field of stochastic processes and their applications. See \cite{sz} and references therein.  Adiabatic times, as introduced in \cite{kov}, is a way to quantify the stability for a certain class of time inhomogeneous Markov processes.  
In order for us to introduce the reader to the type of adiabatic results that we will be working with in this paper, let us first mention earlier results that were published in \cite{kov}, thus postponing a more elaborate discussion of the matter until section \ref{results}.
\vskip 0.2 in
\noindent
Mixing time quantifies the time it takes for a Markov chain to reach a state that is close enough to its stationary distribution.  For the discrete-time finite state case we will look at the evolution of the Markov chain through its probability transition matrix. See \cite{lpw} for a systematized account of mixing time theory and examples. Let $\|\cdot\|_{TV}$ denote the total variation distance.

\begin{Def}
Suppose P is a discrete-time finite Markov chain with a unique stationary distribution $\pi$, i.e. $\pi P = \pi$. Given an $\epsilon > 0$, the \underline{\text{mixing time}} $t_{mix}(\epsilon)$ is defined as $$ t_{mix}(\epsilon) = \inf \left\{ t : \| \nu P^{t} - \pi \|_{TV} \leq \epsilon, \text{ for all probability distributions } \nu \right\}. $$
\end{Def}

\noindent
To define adiabatic time in its first and simplest form (that we will expand and generalize a few pages down) we have to consider a time inhomogeneous Markov chain whose probability transition matrix evolves linearly from an initial probability transition matrix $P_{initial}$ to a final probability transition matrix $P_{final}$.  Namely, we consider two transition probability operators, $P_{initial}$ and $P_{final}$, on a finite state space $\Omega$, and we suppose there is only a unique stationary distribution $\pi_{f}$ of $P_{final}$. 
We let 
\begin{equation} \label{discreteOLD}
P_{s}=(1-s)P_{initial} + sP_{final}
\end{equation}

\noindent
We use (\ref{discreteOLD}) to define a time inhomogeneous Markov chain $P_{t \over T}$ over $[0,T]$ time interval. The adiabatic time quantifies how gradual the transition from $P_{initial}$ to $P_{final}$ should be so that at time $T$, the distribution is $\epsilon$ close to the stationary  distribution $\pi_{f}$ of $P_{final}$. 
\begin{Def}\label{discrete_time}
Given $\epsilon > 0$, a time $T_{\epsilon}$ is called the \underline{\text{adiabatic time}} if it is the least T such that $$max_{\nu} \| \nu P_{\frac{1}{T}}P_{\frac{2}{T}}\cdots P_{\frac{T-1}{T}}P_{1} - \pi_{f} \|_{TV} \leq \epsilon$$
where the maximum is taken over all probability distributions $\nu$ over $\Omega$.
\end{Def}

With these definitions one would naturally ask how adiabatic and mixing times compare.  
This will be especially relevant given the emergence of quantum adiabatic computation and some instances of using adiabatic algorithms to solve certain classical computation problems.  See \cite{siam} and \cite{rss}. It can be speculated that there may be scenarios in which the adiabatic time is more convenient to compute than mixing times.  If we find the relationship between the two, it will give us an understanding of the adiabatic transition (which is more prevalent in a context of physics) in terms of mixing times and vice versa. The following adiabatic theorem was proved in \cite{kov}. 
\begin{thm*}[Kovchegov 2009]
Let $t_{mix}$ denote the mixing time for $P_{final}$. Then the adiabatic time $$T_{\epsilon} = \mathcal{O} \left( \frac{t_{mix}(\epsilon/2)^{2}}{\epsilon} \right)$$
\end{thm*}
\noindent
In subsection \ref{example} we will give an example to show that $t_{mix}^2$ is the best bound for the adiabatic time in this setting. 
There $\Omega=\{0,1,2,\dots,n\}$ and
$$P_{initial} = \left(\begin{array}{cccc}1 & 0 & \cdots & 0 \\1 & 0 & \cdots & 0 \\\vdots & \vdots & \ddots & \vdots \\1 & 0 & \cdots & 0\end{array}\right)
\quad \text{ and } \quad
P_{final}=\left(\begin{array}{cccccc}0 & 1 & 0 & 0 & \cdots & 0 \\0 & 0 & 1 & 0 & \cdots & 0 \\0 & 0 & 0 & 1 & \ddots & \vdots \\\vdots & \vdots & \vdots & \ddots & \ddots & 0 \\0 & 0 & 0 & \cdots & 0 & 1 \\0 & 0 & 0 & \cdots & 0 & 1\end{array}\right)$$
\vskip 0.2 in
\noindent
Similar adiabatic  results hold in the case of continuous-time Markov chains.  There, the concept of an adiabatic time is defined within the same setting and a relationship with mixing time is shown.  Let us state a continuous adiabatic result from \cite{kov}, and then prove a more general statement of the theorem in the next section. 
\vskip 0.2 in
\noindent
Once again we define mixing time as a measurement of the time it takes for a Markov chain to reach a state that is close enough to its stationary distribution.  For the continuous-time, finite-state case we look at the evolution of the Markov chain through its probability transition matrix as a function over time.
\begin{Def}
Suppose $P(t)$ is a finite continuous-time Markov chain with a unique stationary distribution $\pi$. Given an $\epsilon > 0$, the \underline{\text{mixing time}} $t_{mix}(\epsilon)$ is defined as $$t_{mix}(\epsilon) = \inf \left\{ t: \| \nu P(t) - \pi \|_{TV} \leq \epsilon,
\text{ for all probability distributions } \nu \right\}.$$
\end{Def}

\noindent
To define an adiabatic time we have to look at the linear evolution of a generator for the initial probability transition matrix to a generator for the final probability transition matrix.  Suppose $Q_{initial}$ and $Q_{final}$ are two bounded generators for continuous-time Markov processes on a finite state space $\Omega$, and $\pi_{f}$ is the unique stationary distribution for $Q_{final}$. Let us define a time inhomogeneous generator 
\begin{equation} \label{continuousOLD}
Q[s] = (1-s)Q_{initial} + sQ_{final} 
\end{equation}
Given $T > 0$ and $0 \leq t_{1} \leq t_{2} \leq T$, let $P_{T}(t_{1},t_{2})$ denote a matrix of transition probabilities of a Markov process generated by $Q[\frac{t}{T}]$ over the time interval $[t_{1},t_{2}]$. 
With this new generator we define the adiabatic time to be the smallest transition time $T$ such that regardless of our starting distribution, the continuous-time Markov chain generated by  $Q[\frac{t}{T}]$ arrives at a state close enough to our stationary distribution $\pi_{f}$.

\begin{Def}\label{continuous_time}
Given $\epsilon >0$, a time $T_{\epsilon}$ is called the \underline{\text{adiabatic time}} if it is the least T such that $$ \max_{\nu} \| \nu P_{T}(0,T) - \pi_{f} \|_{TV} \leq \epsilon$$
where the maximum is taken over all probability distributions $\nu$ over $\Omega$.
\end{Def}
The above definition for continuous-time Markov chains is similar to the one in the discrete time setting. The corresponding adiabatic theorem for the continuous times case was proved in \cite{kov}.

\begin{thm*}[Kovchegov 2009]
Let $t_{mix}$ denote the mixing time for $Q_{final}$. Take $\lambda$ such that $\lambda \geq \max_{i \in \Omega} \sum_{j:j \neq i}q_{i,j}^{initial}$
 and  $ \lambda \geq \max_{i \in \Omega} \sum_{j:j \neq i}q_{i,j}^{final}$,
where $q_{i,j}^{initial}$ and $q_{i,j}^{final}$ are the rates in $Q_{initial}$ and $Q_{final}$ respectively. Then the adiabatic time $$T_{\epsilon} \leq \frac{\lambda t_{mix}(\epsilon/2)^{2}}{\epsilon} + \theta$$
where $\theta = t_{mix}(\epsilon/2) + \epsilon/(4\lambda)$.
\end{thm*}
\noindent
This is once again the best bound as can be shown through the corresponding example. 
\vskip 0.2 in
\noindent
In the next section we will state the adiabatic results for Markov chains that generalize the above mentioned theorems in \cite{kov} and provide examples of applications in statistical mechanics. Section \ref{proofs} is dedicated to proofs.

\section{Results and applications}\label{results}

Here we extend the results from \cite{kov}, and thus expand the range of problems that can be analyzed with these types of adiabatic theorems.  One such problem that we will discuss in subsection \ref{Ising} deals with adiabatic Glauber dynamics for Ising model.  Now, in order to solve a larger class of problems, we redefine the adiabatic transition for both the discrete and continuous cases. \\

\noindent
We consider an adiabatic dynamics where transition probabilities change gradually from $P_{initial}=\left\{p_{i,j}^{initial}\right\}$ to $P_{final}=\left\{p_{i,j}^{final}\right\}$ so that, for each pair of states $i$ and $j$, the corresponding mutation of $p_{i,j}$ from $p_{i,j}^{initial}$ to $p_{i,j}^{final}$ is implemented differently and not always linearly. 
In the case of discrete time steps, this means defining
\begin{equation}\label{discreteNEW}
p_{i,j}[s] = (1-\phi_{i,j}(s)) p_{i,j}^{initial} + \phi_{i,j}(s) p_{i,j}^{final}~,
\end{equation}
where $\phi_{i,j}:[0,1] \rightarrow [0,1]$ are continuous functions such that $\phi_{i,j}(0)=0$ and $\phi_{i,j}(1)=1$  for all locations (i,j).
\vskip 0.2 in
\noindent
The above definition generalizes (\ref{discreteOLD}). If we suppose there is a unique stationary distribution $\pi_{f}$ for $P_{final}$, then the Definition \ref{discrete_time} of {\bf adiabatic time}  $T_{\epsilon}$ given in the previous section will hold for the adiabatic dynamics defined in (\ref{discreteNEW}). The new $T_{\epsilon}$ is related to mixing time via the following adiabatic theorem, that we will prove in section \ref{proofs}.

\begin{thm}[Discrete Adiabatic Theorem] \label{discrete}
Let $P_{t \over T}=\left\{p_{i,j}\left[{t \over T}\right]\right\}$ be an inhomogeneous discrete-time Markov chain over $[0,T]$. 
Let $\phi(s)=\min_{i,j} \phi_{i,j}(s)$  be the pointwise minimum function of all of the $\phi_{i,j}$  functions. 
If $m \geq 1$ is an integer such that $\phi$ is $m+1$ times continuously differentiable in a neighborhood of $1$,
$$\phi^{(k)}(1)=0\text{ for all integers }k\text{ such that }1 \leq k <m$$
 and $\phi^{(m)}(1) \neq 0$, then 
$$T_{\epsilon} = \mathcal{O} \left( \frac{t_{mix}^{\frac{m+1}{m}}(\epsilon/2)}{\epsilon^{\frac{1}{m}}} \right)$$
\end{thm}

\vskip 0.2 in
\noindent
The above is, in fact, the best bound in the new setting as shown through the example given later. See subsection \ref{example}.

\vskip 0.2 in
\noindent
Now we extend the notion of adiabatic dynamics for the continuous-time Markov generators as follows. We let
\begin{equation}\label{continuousNEW}
q_{i,j}[s] = (1-\phi_{i,j}(s)) q_{i,j}^{initial} + \phi_{i,j}(s) q_{i,j}^{final} \quad \text{ for all pairs }~i \not=j,
\end{equation}
where once again $\phi_{i,j}:[0,1] \rightarrow [0,1]$ are continuous functions such that $\phi_{i,j}(0)=0$ and $\phi_{i,j}(1)=1$  for all locations (i,j). Also, we let $Q[s]$ denote the corresponding Markov operator. 
\vskip 0.2 in
\noindent
If there is a unique stationary distribution $\pi_{f}$ for $Q_{final}$, then the Definition \ref{continuous_time} of {\bf adiabatic time}  will apply for the extended adiabatic dynamics in (\ref{continuousNEW}), and the new $T_{\epsilon}$ can be  again related to mixing time.
\begin{thm}[Continuous Adiabatic Theorem] \label{continuous}
Let $Q\left[{t \over T}\right]~(t \in[0,T])$ generate the inhomogeneous discrete-time Markov chain. 
Let $\phi(s)=\min_{i,j} \phi_{i,j}(s)$  be the pointwise minimum function of all of the $\phi_{i,j}$  functions. Suppose $m \geq 1$ is an integer such that $\phi$ is $m+1$ times continuously differentiable in a neighborhood of $1$,
$$\phi^{(k)}(1)=0\text{ for all integers }k\text{ such that }1 \leq k <m$$
 and $\phi^{(m)}(1) \neq 0$. If we take $\lambda$ such that 
$$\lambda \geq \max_{i \in \Omega} \sum_{j:j \neq i}q_{i,j}^{initial}~~~~\text{ and  }~~~~\lambda \geq \max_{i \in \Omega} \sum_{j:j \neq i}q_{i,j}^{final},$$
where $q_{i,j}^{initial}$ and $q_{i,j}^{final}$ are the rates in $Q_{initial}$ and $Q_{final}$ respectively. 
Then $$T_{\epsilon} = \mathcal{O} \left( \left[{\lambda \over \epsilon}\right]^{1 \over m} t_{mix}^{\frac{m+1}{m}}(\epsilon/2) \right)$$
\end{thm}
\noindent
The reader can reference the proof of this theorem in section \ref{proofs}.  Again this is the best bound in the new setting as can be shown through the same example.  See subsection \ref{example}.
\vskip 0.2 in
\noindent
Now we check that the above continuous adiabatic theorem is scale invariant. For a positive $M$, we scale the initial and final generators to be ${1 \over M}Q_{initial}$ and ${1 \over M} Q_{final}$ respectively. Then the adiabatic evolution is slowed down $M$ times, and the new adiabatic time should be of order $M\left[{\lambda \over \epsilon}\right]^{1 \over m} t_{mix}^{\frac{m+1}{m}}(\epsilon/2)$ with the old $t_{mix}$ and $\lambda$ taken before scaling. On the other hand the new mixing time will be $Mt_{mix}$, and the new $\lambda$ is ${\lambda \over M}$ as the rates are $M$ times lower. Plugging the new parameters into the expression in the theorem, we obtain
$$\left[{\lambda \over M\epsilon}\right]^{1 \over m} (Mt_{mix})^{\frac{m+1}{m}}=M\left[{\lambda \over \epsilon}\right]^{1 \over m} t_{mix}^{\frac{m+1}{m}}$$
confirming the theorem is invariant under time scaling.
\vskip 0.2 in
\noindent
Let us revisit adiabatic theorems in physics and quantum mechanics. The reader can find a version of quantum adiabatic theorem in \cite{m} and multiple other sources. 
\vskip 0.2 in
\noindent
The adiabatic results in physics consider a system that transitions from one state to another, while the energy function changes from an initial ${\cal H}_{initial}$ to ${\cal H}_{final}$. If the change in the energy function happens slowly enough, for the system that is initially at one of the equilibrium states (i.e. at an eigenvector of the initial energy function ${\cal H}_{initial}$), the resulting state of the transition of the system will end up at a state that is $\epsilon$ close to the corresponding eigenvector of the final energy function ${\cal H}_{final}$. 
That is, provided the change in the external conditions is gradual enough, the $j$th eigenstate of ${\cal H}_{initial}$ is carried to an $\epsilon$ proximity of the $j$th eigenstate of ${\cal H}_{final}$.
\vskip 0.2 in
\noindent
Often the adiabatic results concern with one eigenstate, the ground state. Thinking of Schr\"odinger equation as an $\ell^2(\mathbb{C}^n)$ norm preserving linear dynamics, and a finite Markov process as a natural description of an $\ell^1(\mathbb{R}^n_+)$ norm preserving linear dynamics,  the ground state of one would correspond to the stationary state of the other.  It is important to mention that in addition to all above properties, the quantum adiabatic theorems often require the transition to be gradual enough for the state to be within an $\epsilon$ proximity of the corresponding ground state at each time during the transition. Taking this into account, the complete analogue of quantum adiabatic theorem for  $\ell^1(\mathbb{R}^n_+)$ would be the one in which 
the initial distribution is $\mu_0=\pi_{initial}$ and 
$$\| \mu_t - \pi_t \| < \epsilon~~~\forall t \in [0,T],$$
where $\mu_t=\mu_0P_{\frac{1}{T}}P_{\frac{2}{T}}\cdots P_{\frac{T-1}{T}}P_{1}$ is the distribution of the inhomogeneous Markov chain at time $t \in [0,T]$, $\pi_{initial}$ is the stationary distribution of $P_{initial}$, and $\pi_t$ is  the stationary distribution $P_{t \over T}$. While we are currently working on proving the above mentioned complete analogue in both discrete and continuous cases, the adiabatic results of this section are sufficiently strong for answering our questions concerning adiabatic Glauber dynamics as stated in the following subsection.

\subsection{Applications to Ising models with adiabatic Glauber dynamics}\label{Ising}

Let us first state a version of the quantum  adiabatic theorem.
 Given two Hamiltonians, ${\cal H}_{initial}$ and ${\cal H}_{final}$, acting on a quantum system. 
 Let 
 \begin{equation} \label{hamiltonian}
 {\cal H}(s)=(1-s){\cal H}_{initial}+s{\cal H}_{final}
 \end{equation}
 Suppose the system evolves according to ${\cal H}(t/T)$ from time $t=0$ to time $T$ .
Then if $T$ is large enough, the final state of the system will be close to
the ground state of ${\cal H}_{final}$. They are $\epsilon$ close in the $\ell^2$ norm whenever $T \geq {C \over \epsilon \beta^3}$,
where $\beta$ is the least spectral gap of ${\cal H}(s)$ over all $s \in [0,1]$, and $C$ depends linearly on a square of the distance between ${\cal H}_{initial}$ and ${\cal H}_{final}$.
\vskip 0.2 in
\noindent
Now, switching to canonical ensembles of statistical mechanics will land us in a Gibbs measure space with familiar probabilistic properties, i.e.  the Markov property of statistical independence.
We consider a nearest-neighbor Ising model.  There the spins can be of two types, -1 and +1.  The spins interact only with nearest neighbors.  A Hamiltonian determines the energy-value of the interactions of the configuration of spins. 
\vskip 0.2 in
\noindent
Here, for a microstate, we multiply its energy by the thermodynamic beta and call it the Hamiltonian of the microstate. In other words, letting $x$  be a configuration of spins, the Hamiltonian we use in this paper will be defined as 
$$\mathcal{H}(x) = -{\beta \over 2} \sum_{i \not= j}  M_{i,j} x(i) x(j)$$ 

\noindent
where $\beta$  is the thermodynamic beta, i.e. its inverse is the temperature times  Boltzmann's constant, $M=\{M_{i,j}\}$  is a symmetric matrix and for locations $i$ and $j$, $M_{i,j} = 0$  if $i$ is not a nearest neighbor to $j$ and $M_{i,j}=1$  if $i$ is a  nearest neighbor to $j$. 
\vskip 0.2 in
\noindent
The Markov property of statistical independence is reflected through the {\it local Hamiltonian} defined at every location $j$ as follows
$$ \mathcal{H}^{loc}(x(j)) = -\beta \sum_{i:i \sim j} x(i)x(j),$$
where $i \sim j$ means $i$ and $j$ are nearest neighbors on the graph.
\vskip 0.2 in
\noindent
 In the original, non-adiabatic case, the Glauber dynamics is used to generate the following Gibbs distribution
 $$\pi(x)={1 \over Z(\beta)}e^{-{\cal H}(x)}$$
 over all spin configurations $x \in \{-1,+1\}^S$, 
 where $S$ denotes all the sites of a graph,  and $Z(\beta)$ is the normalization constant.
 Let us describe how the Glauber dynamics works in the case when each vertex of the connected graph is of the same degree. There for each location $j$ we have an independent exponential clock with parameter one associated with it. When the clock rings, the spin $x(j)$ of configuration $x$ at the site $j$ on the graph is reselected using the following probability 
 $$P(x(j)=+1)=\frac{e^{-\mathcal{H}^{loc}(x_{+}(j))}}{e^{-\mathcal{H}^{loc}(x_{-}(j))}+e^{-\mathcal{H}^{loc}(x_{+}(j))}}=2-2\tanh\left\{\mathcal{H}^{loc}(x_{+}(j))\right\}$$
where $x_{+}(i) = x_{-}(i) = x(i)$ for $i \neq j$, $x_{+}(j)=+1$ and $x_{-}(j)=-1$. Here $$P(x(j)=-1)=1-P(x(j)=+1)$$ Also 
$\mathcal{H}^{loc}(x_{-}(j)) = -\mathcal{H}^{loc}(x_{+}(j))$. Now we have a continuous-time Markov process, in which our states are the collection of the configurations of spins.  
\vskip 0.2 in
\noindent
Now, consider an adiabatic evolution of Hamiltonians as in (\ref{hamiltonian}). There at each time $t$, $${\cal H}(s)=(1-s){\cal H}_{initial}+s{\cal H}_{final},$$
where $s={t \over T}$. The local Hamiltonians must therefore evolve accordingly,
$${\cal H}^{loc}_{s}=(1-s){\cal H}_{initial}^{loc}+s{\cal H}_{final}^{loc}$$
and the {\it adiabatic Glauber dynamics} is the one where when the clock rings, the spin $x(j)$ is reselected with probabilities 
$$P_{s}(x(j)=+1)=\frac{e^{-\mathcal{H}_{s}^{loc}(x_{+}(j))}}{e^{-\mathcal{H}_{s}^{loc}(x_{-}(j))}+e^{-\mathcal{H}_{s}^{loc}(x_{+}(j))}}$$
Here too, $\mathcal{H}_{s}^{loc}(x_{-}(j)) = -\mathcal{H}_{s}^{loc}(x_{+}(j))$. 
\vskip 0.2 in
\noindent
The stationary distribution of the $Q_{final}$-generated Markov process, i.e. Glauber dynamics with ${\cal H}_{initial}$ energy,  is, for a configuration $x$, 
$$ \pi(x) = \frac{e^{- \mathcal{H}_{final}(x)}}{\sum_{\text{all config. }x'}e^{- \mathcal{H}_{final}(x')}} $$

\subsubsection{Adiabatic Glauber dynamics on $\mathbb{Z}^2 /n \mathbb{Z}^2$}
Consider nonlinear adiabatic Glauber dynamics of an Ising model on a two-dimsnsional torus $\mathbb{Z}^2 /n \mathbb{Z}^2$.
There any two neighboring spin configurations $x$ and $y$ in $\{-1,+1\}^{n^2}$ differ at only one site on the graph, say $v \in \mathbb{Z}^2 /n \mathbb{Z}^2$. That is $y(u)=\begin{cases}
     x(u) & \text{ if } u \not= v \\
     1-x(u) & \text{ if }u=v
\end{cases}$.
The transition rates evolve according to the adiabatic Glauber dynamics rules, and the transition rates can be represented as
$$q_{x,y}[s] = (1-\phi_{x,y}(s)) q_{x,y}^{initial} + \phi_{x,y}(s) q_{x,y}^{final}$$
as in (\ref{continuousNEW}). Here the functions $\phi_{x,y}(s)$ for two neighbors $x$ and $y$ depend entirely on the spins around the discrepancy site $v$. Namely if all four neighbors of $v$ are of the same spin ($+1$ or $-1$), then
$$\phi_{x,y}(s)= \frac{cosh(-4 \beta_{2}) \cdot sinh(s(4 \beta_{1}-4 \beta_{2}))}{sinh(4 \beta_{1} - 4 \beta_{2}) \cdot cosh(-4 \beta_{1} + s(4 \beta_{1}- 4 \beta_{2}))}$$
If it is three of one kind, and one of the other (i.e three $+1$ and one $-1$, or three $-1$ and one $+1$) as illustrated below
$$\left.\begin{array}{ccccc}  &   & -1 &   &   \\  &   & | &   &   \\+1 & - & v & - & +1  \\  &   & | &   &   \\  &   & +1 &   &  \end{array}\right.$$
then
$$\phi_{x,y}(s)=\frac{cosh(-2 \beta_{2}) \cdot sinh(s(2 \beta_{1}-2 \beta_{2}))}{sinh(2 \beta_{1} - 2 \beta_{2}) \cdot cosh(-2 \beta_{1} + s(2 \beta_{1}-2 \beta_{2}))} $$
\vskip 0.2 in
\noindent
If there are two of each kind, any function works, as both, the initial and the final, local Hamiltonians produce the same transition rates $q_{x,y}^{initial} = q_{x,y}^{final}=1/2$.
\vskip 0.2 in
\noindent
Now, since
$$\frac{cosh(-4 \beta_{2}) \cdot sinh(s(4 \beta_{1}-4 \beta_{2}))}{sinh(4 \beta_{1} - 4 \beta_{2}) \cdot cosh(-4 \beta_{1} + s(4 \beta_{1}- 4 \beta_{2}))} \geq \frac{cosh(-2 \beta_{2}) \cdot sinh(s(2 \beta_{1}-2 \beta_{2}))}{sinh(2 \beta_{1} - 2 \beta_{2}) \cdot cosh(-2 \beta_{1} + s(2 \beta_{1}-2 \beta_{2}))}$$
for $s \in [0,1]$, Theorem \ref{continuous} implies the adiabatic time
$$ T_{\epsilon} = \mathcal{O}\left(C\frac{n^2}{\epsilon}\left[\log(n) + \log\left(\frac{2}{\epsilon}\right)\right]^{2}\right),$$
where $C={(2 \beta_{1} - 2 \beta_{2})\left[coth(2 \beta_{1} - 2 \beta_{2})-tanh(-2 \beta_{2})\right] \over [1-tanh(2\beta_{2})]^2}$.
Here, at every vertex on the torus we attached a Poisson clock with rate  one, and therefore we can take $\lambda=n^2$.
Also $m=1$ in the theorem, and one can find the expression for $t_{mix}$ in \cite{lpw}.

\subsubsection{Adiabatic Glauber dynamics on $\mathbb{Z}^d /n \mathbb{Z}^d$}
The adiabatic Glauber dynamics of an Ising model on a $d$-dimsnsional torus $\mathbb{Z}^d /n \mathbb{Z}^d$ solves similarly.
There the minimum function $\phi(s)$ of the Theorem \ref{continuous} is same as in the case of $d=2$
$$\phi(s) = \frac{cosh(-2 \beta_{2}) \cdot sinh(s(2 \beta_{1}-2 \beta_{2}))}{sinh(2 \beta_{1} - 2 \beta_{2}) \cdot cosh(-2 \beta_{1} + s(2 \beta_{1}-2 \beta_{2}))}$$
and the adiabatic time
$$T_{\epsilon} = \mathcal{O}\left(C\frac{n^d}{\epsilon}\left[\log(n) + \log\left(\frac{2}{\epsilon}\right)\right]^{2}\right),$$
where again $C={(2 \beta_{1} - 2 \beta_{2})\left[coth(2 \beta_{1} - 2 \beta_{2})-tanh(-2 \beta_{2})\right] \over [1-tanh(2\beta_{2})]^2}$.
\vskip 0.2 in
\noindent
Notice that the time scaling argument that followed the statement of Theorem \ref{continuous} works here as well. That is, if we use one Poisson clock of rate one for all vertices, or equivalently place Poisson clocks of rate $n^{-d}$ at every individual vertex, the new adiabatic time will be 
$$T'_{\epsilon} = \mathcal{O}\left(C\frac{n^{2d}}{\epsilon}\left[\log(n) + \log\left(\frac{2}{\epsilon}\right)\right]^{2}\right)$$
as $\lambda'=1$ here.

\subsection{The bound is optimal}\label{example}
In this subsection we give examples that show that the $t_{mix}^2$ order of adiabatic time given in Kovchegov \cite{kov}, and $t_{mix}^{m+1 \over m}$ order for more general settings of this current paper are in fact optimal. We consider discrete parobability transition matrices
$$P_{initial} = \left(\begin{array}{cccc}1 & 0 & \cdots & 0 \\1 & 0 & \cdots & 0 \\\vdots & \vdots & \ddots & \vdots \\1 & 0 & \cdots & 0\end{array}\right)
\quad \text{ and } \quad
P_{final}=\left(\begin{array}{cccccc}0 & 1 & 0 & 0 & \cdots & 0 \\0 & 0 & 1 & 0 & \cdots & 0 \\0 & 0 & 0 & 1 & \ddots & \vdots \\\vdots & \vdots & \vdots & \ddots & \ddots & 0 \\0 & 0 & 0 & \cdots & 0 & 1 \\0 & 0 & 0 & \cdots & 0 & 1\end{array}\right)$$
over $n+1$ states, $\{0,1,2,\dots,n\}$, and let the  discrete-time adiabatic probability transition matrix to be
$$P_{s}=(1-s)P_{initial} + sP_{final}$$ 
as in \cite{kov}. Let $\pi_f$ again denote the stationary distribution for $P_{final}$. Here $\pi_f=(0,\dots,0,1)$ and the mixing time $t_{mix}(\epsilon)=n$ for any $\epsilon \in (0,1)$. 
\vskip 0.2 in
\noindent
Now, since $\mu P_{initial}=\rho$ for any probability distribution $\mu$,  we have the following inequality

$$\| \nu P_{\frac{1}{T}} \cdots P_{\frac{T}{T}} - \pi_{f} \|_{TV} \geq \left\| \rho \left( \sum_{j=0}^{T} \left( 1 - \frac{j}{T} \right) \frac{T!}{j! \cdot T^{T-j}} \right) \left( P_{fin}^{T-j}  - \pi_{f} \right) \right\|_{TV}$$ 

\noindent
Observe that $\rho P_{fin}^{T-j} = \pi_{f}$ for any $0 \leq j \leq T-n$. Therefore
\begin{align*}
\| \nu P_{\frac{1}{T}} \cdots P_{\frac{T}{T}} - \pi_{f} \|_{TV} 
&\geq \sum_{j=T-n+1}^{T} \left( 1- \frac{j}{T} \right) \frac{T!}{j! T^{T-j}} \\
&= \sum_{j=T-n+1}^{T} \left( \frac{T!}{j! \cdot T^{T-j}} - \frac{T!}{(j-1)! \cdot T^{T-(j-1)}} \right) \\
&= 1 - \frac{T!}{(T-n)! \cdot T^{n}} \\
&=1 - \frac{T-n+1}{T} \cdots \frac{T-1}{T}.
\end{align*}

\noindent
Now, because $\frac{T-n+1}{T} \cdots \frac{T-1}{T} \leq \left( \frac{T - \frac{n}{2}}{T} \right)^{\frac{n}{2}}$  for $n \geq 2$, we see that 
$$\| \nu P_{\frac{1}{T}} \cdots P_{\frac{T}{T}} - \pi_{f} \|_{TV} \geq 1 - \left( \frac{T - \frac{n}{2}}{T} \right)^{\frac{n}{2}} 
\geq1 - e^{-\left( \frac{n^{2}}{4T} \right)}$$

\noindent
Thus $\epsilon \geq \| \nu P_{\frac{1}{T}} \cdots P_{\frac{T}{T}} - \pi_{f} \|_{TV}  \geq 1 - e^{-\left( \frac{n^{2}}{4T} \right)}$ implies
$T \geq  \frac{-n^{2}}{4\log(1- \epsilon)} \approx  \frac{n^{2}}{4 \epsilon}={t_{mix}^2 \over 4\epsilon} $, proving that the order of adiabatic time $T_{\epsilon} = \mathcal{O} \left( \frac{t_{mix}(\epsilon/2)^{2}}{\epsilon} \right)$ in \cite{kov} is optimal.

 \subsubsection{Optimal bound for general settings}
The same example works in the more general setting considered in this paper.  For the same $P_{initial}$ and $P_{final}$, let
$$p_{i,j}[s] = (1-\phi_{i,j}(s)) p_{i,j}^{initial} + \phi_{i,j}(s) p_{i,j}^{final}$$
as in (\ref{continuousNEW}). 
Suppose $\phi_{i,j}(s)= \phi(s)$  for all pairs of states $(i,j)$, and suppose $m \geq 1$ is an integer such that $\phi$ is $m+1$ times continuously differentiable in a neighborhood of $1$,
$$\phi^{(k)}(1)=0\text{ for all integers }k\text{ such that }1 \leq k <m$$
 and $\phi^{(m)}(1) \neq 0$.
Then
 $$ \| \nu P_{\frac{1}{T}} \cdots P_{\frac{T}{T}} - \pi_{f} \|_{TV} \geq \left\| \rho \left( \sum_{l=0}^{T-1} \left( 1 - \phi(l/T) \right) \prod_{j=l+1}^{T}\phi(j/T) \right) \left( P_{final}^{T-l}  - \pi_{f} \right) \right\|_{TV}$$
 and therefore
 \begin{align*}
\| \nu P_{\frac{1}{T}} \cdots P_{\frac{T}{T}} - \pi_{f} \|_{TV} 
&\geq \sum_{l=T-n+1}^{T-1} \left( 1- \phi(l/T) \right) \prod_{j=l+1}^{T} \phi(j/T) \\
&= \sum_{l=T-n+1}^{T-1} \left( \prod_{j=l+1}^{T}\phi(j/T) - \prod_{j=l}^{T}\phi(j/T) \right) \\
&= 1 - \prod_{j=T-n+1}^{T} \phi(j/T). 
\end{align*}

\noindent
The minimum function $\phi(x) = 1 +\frac{\phi^{(m)}(1)(x-1)^{m}}{m!}+O\left(|x-1|^{m+1}\right)$ and
\begin{align*}
\| \nu P_{\frac{1}{T}} \cdots P_{\frac{T}{T}} - \pi_{f} \|_{TV} 
&\geq  1 - \prod_{j=T-n+1}^{T} \left(1 + \frac{(-1)^{m}\phi^{(m)}(1) \cdot (T-j)^{m}}{T^{m} \cdot m!}+O\left((1-j/T)^{m+1}\right) \right) \\
&= 1 - e^{\sum_{j=T-n+1}^{T} \log \left( 1 + \frac{(-1)^{m}\phi^{(m)}(1) \cdot (T-j)^{m}}{T^{m} \cdot m!}+O\left((1-j/T)^{m+1}\right) \right)}\\
&\geq 1-e^{{(-1)^{m}\phi^{(m)}(1) \over T^m}\cdot \sum_{j=1}^{n-1}j^m+O((n/T)^{m+1})}
\end{align*}
as $\log(1+x) \leq x$.
\vskip 0.2 in
\noindent
It is a well know fact that 
\begin{equation} \label{Bernoulli}
\sum_{j=1}^{n-1}j^{k} = \sum_{j=0}^{k}{\frac{B_{j}}{(k+1)-j}\binom{k}{j} n^{(k+1)-j}},
\end{equation}
where $B_{j}$ is the $j$th Bernoulli number. 
\vskip 0.2 in
\noindent
Suppose $\epsilon \geq \| \nu P_{\frac{1}{T}} \cdots P_{\frac{T}{T}} - \pi_{f} \|_{TV}$, then 
$$ \epsilon \approx -\log (1- \epsilon) \geq {(-1)^{m+1}\phi^{(m)}(1) \over T^m}\cdot \sum_{j=0}^{m}{\frac{B_{j}}{(m+1)-j} \binom{m}{j} n^{(m+1)-j}}+O((n/T)^{m+1})$$
Thus confirming that the order of adiabatic time $T_{\epsilon} = \mathcal{O} \left( \frac{t_{mix}^{\frac{m+1}{m}}}{\epsilon^{\frac{1}{m}}} \right)$ in Theorem \ref{discrete} is optimal.
\vskip 0.2 in
\noindent
Naturally, there is a similar example in the continuous case. There
$$Q_{initial} = \left(\begin{array}{ccccc}0 & 0 & 0 & \cdots & 0 \\1 & -1 & 0 & \cdots & 0 \\1 & 0 & -1 & \ddots & \vdots \\\vdots & \vdots & \ddots & \ddots & 0 \\1 & 0 & \cdots & 0 & -1\end{array}\right)
\quad \text{ and } \quad
Q_{final}=\left(\begin{array}{cccccc}-1 & 1 & 0 & 0 & \cdots & 0 \\0 & -1 & 1 & 0 & \cdots & 0 \\0 & 0 & -1 & 1 & \ddots & \vdots \\\vdots & \vdots & \ddots & \ddots & \ddots & 0 \\0 & 0 & 0 & \cdots & -1 & 1 \\0 & 0 & 0 & \cdots & 0 & 0\end{array}\right)$$

 \section{Proofs}\label{proofs}
In this section we give formal proofs to both adiabatic theorems, Theorem \ref{discrete} and Theorem \ref{continuous}.

\subsection{Proof of Theorem \ref{discrete}}

\begin{proof}
We write 
$$ p_{i,j}[s] = (1-\phi_{i,j}(s)) p_{i,j}^{(initial)} + (\phi_{i,j}(s) - \phi(s)) p_{i,j}^{(final)} + \phi(s) p_{i,j}^{(final)}$$
and define transition probability matrix $\hat{P}$ to be  such that 
$$ (1- \phi(s)) \hat{p}_{i,j} = (1-\phi_{i,j}(s)) p_{i,j}^{(initial)} + (\phi_{i,j}(s) - \phi(s)) p_{i,j}^{(final)}$$ 

\noindent
We will thus have that $$ P_{s} = (1-\phi(s)) \hat{P} + \phi(s) P_{final}. $$ 

\noindent
Observe that 
$$\nu P_{\frac{1}{T}} P_{\frac{2}{T}} \cdots P_{\frac{T-1}{T}} P_{1} = \left[ \prod_{j=N+1}^{T} \phi(j/T)\right]\nu_{N} P_{final}^{T-N} + \cal{E}$$
where $\nu_{N} = \nu P_{\frac{1}{T}}P_{\frac{2}{T}} \cdots P_{\frac{N}{T}}$, and $\cal{E}$ is the rest of the terms, and both $T$ and $N$ are natural numbers with $N < T$.  

\noindent
By the triangle inequality, we have $$ \max_{\nu}{\| \nu P_{\frac{1}{T}} P_{\frac{2}{T}} \cdots P_{\frac{T-1}{T}} P_{1} - \pi_{f} \|_{TV}} \leq \max_{\nu}{\| \nu P_{final}^{T-N} - \pi_{f} \|_{TV}} \cdot \left[ \prod_{j=N+1}^{T} \phi(j/T) \right] + S_{N} $$
where $0 \leq S_{N} \leq 1 - \left[ \prod_{j=N+1}^{T} \phi(\frac{j}{T}) \right]$. \\

\noindent
Let set $T-N = t_{mix}(\epsilon/2)$, where $\epsilon > 0$  is small. Then we have that 
$$ \max_{\nu}{\| \nu P_{final}^{T-N} - \pi_{f} \|_{TV}} \cdot \left[ \prod_{j=N+1}^{T} \phi(j/T) \right] \leq \epsilon/2$$
 
\noindent
Setting $1 - \left[ \prod_{j=N+1}^{T} \phi(\frac{j}{T}) \right] \leq \epsilon/2$ we obtain
$$ \log{(1 - \epsilon/2)} \leq \sum_{j=N+1}^{T} \log{\phi(j/T)}$$

\noindent
We plug in the approximation of the minimum function $\phi$ around $x=1$
$$\phi(x) = 1 +\frac{\phi^{(m)}(1)(x-1)^{m}}{m!}+O\left(|x-1|^{m+1}\right)$$
obtaining 
$$-\log{(1 - \epsilon/2)} \geq -\sum_{j=N+1}^{T} \log{\left( 1 +\frac{(-1)^m\phi^{(m)}(1)(T-j)^{m}}{T^m \cdot m!}+O\left((1-j/T)^{m+1}\right) \right)}$$
Therefore
$$-\log{(1 - \epsilon/2)} \geq \frac{(-1)^{m+1}\phi^{(m)}(1)}{T^m \cdot m!}\sum_{j=1}^{T-N-1}j^m+O\left({(T-N)^{m+2} \over T^{m+1}}\right)$$
Observe that $(-1)^{m+1}\phi^{(m)}(1) \geq 0$ as $~\phi:[0,1] \rightarrow [0,1]$ and $\phi(1)=1$. 

\noindent
By (\ref{Bernoulli}), 
$~~\sum_{j=1}^{t_{mix}(\epsilon/2)-1}j^{m} = \sum_{k=0}^{m}{\frac{B_{k}}{(m+1)-k} \binom{m}{k} t_{mix}(\epsilon/2)^{(m+1)-k}},~~$
where $B_{k}$ is the $k$th Bernoulli number, and therefore
$$\epsilon >-\log{(1 - \epsilon/2)} \geq  \frac{(-1)^{m+1}\phi^{(m)}(1)}{T^m \cdot m!}\sum_{k=0}^{m}{\frac{B_{k}}{(m+1)-k} \binom{m}{k} t_{mix}(\epsilon/2)^{(m+1)-k}}+O\left({(T-N)^{m+2} \over T^{m+1}}\right)$$
In order for the right hand side of the above equation to be $-\log{(1 - \epsilon/2)}$ close to zero, it is sufficient for $T$ to be of order of 
$\mathcal{O} \left( \frac{t_{mix}^{\frac{m+1}{m}}(\epsilon/2)}{\epsilon^{\frac{1}{m}}} \right)$.
\end{proof}

\subsection{Proof of Theorem \ref{continuous}}

\begin{proof}
Define $\hat{Q}$ to be a Markov generate with off-diagonal entries
$$\hat{q}_{i,j} = {1-\phi_{i,j}(s) \over 1- \phi(s)} q_{i,j}^{(initial)} + {\phi_{i,j}(s) - \phi(s) \over 1- \phi(s)} q_{i,j}^{(final)}$$
Then writing
$$ q_{i,j}[s] = (1-\phi_{i,j}(s)) q_{i,j}^{(initial)} + (\phi_{i,j}(s) - \phi(s)) q_{i,j}^{(final)} + \phi(s) q_{i,j}^{(final)}$$
would imply 
$$ Q[s] = (1-\phi(s)) \hat{Q} + \phi(s) Q_{final}$$
Observe  that 
$$ \lambda \geq \max_{i \in \Omega}{\sum_{j:j \neq i} \hat{q}_{i,j}}~~\text{ and }~~ \lambda \geq \max_{i \in \Omega}{\sum_{j:j \neq i} q_{i,j}\left[ \frac{t}{T} \right] } $$
as 
$$ \lambda \geq \max_{i \in \Omega}{\sum_{j:j \neq i} q_{i,j}^{(initial)}} ~~\text{ and }~~ \lambda \geq \max_{i \in \Omega}{\sum_{j:j \neq i} q_{i,j}^{(final)}}$$

\noindent
Let $P_{final}(t) = e^{tQ_{final}}$ denote the transition probability matrix associated with the generator $Q_{final}$, and let 
$P_{0} = I + \frac{1}{\lambda}\hat{Q}$ and $P_{1} = I + \frac{1}{\lambda}Q_{final}$. 

\noindent
The $P_{0}$ and $P_{1}$ are discrete Markov chains. Conditioning on the number of arrivals within the $[N,T]$ time interval
$$ \nu P_{T}(0,T) = \nu_{N} P_{T}(N,T) = \nu_{N} \left( \sum_{n=0}^{\infty} \frac{(\lambda(T-N))^{n}}{n!}e^{-\lambda (T-N)}I_{n} \right) $$
where $\nu_{N} = \nu P_{T}(0,N)$ and 
\begin{align*}
I_{n} =\frac{n!}{(T-N)^{n}}~~~ \int\!\!\cdots\!\!\int_{\substack{
            ~\\ ~\\ ~\\ \!\!\!\!\!\!\!\!\!\!\!\!\!\!\!\!\!\!\!\!\!\!\!\!\!\!\!\!\!\!\!\ N<s_{1}< \cdots < s_{n} < T}} \!\!\!\!\!\!\!\!\! & \left[ \left(1-\phi\left(\frac{s_{1}}{T}\right)\right)P_{0} + \phi\left(\frac{s_{1}}{T}\right) P_{1}\right]\\
& \qquad \cdots  \left[ \left(1-\phi\left(\frac{s_{n}}{T}\right)\right)P_{0} + \phi\left(\frac{s_{n}}{T}\right) P_{1} \right] ds_{1} \cdots ds_{n}
\end{align*}
\noindent
i.e. the order statistics of the $n$. 


\noindent
Therefore, combining the terms with $P_{final}$, we obtain
\begin{align*} 
\nu P_{T}(0,T) & = \nu_{N} \left( \sum_{n=0}^{\infty} \frac{ \lambda^n P_{final}^{n}}{n!}e^{-\lambda(T-N)} \int_{N}^{T}\!\!\!\!\cdots\!\!\int_{N}^{T} \phi\left(\frac{s_{1}}{T}\right) \cdots \phi\left(\frac{s_{n}}{T}\right) ds_{1} \cdots ds_{n} \right) + \mathcal{E}\\
&= e^{-\lambda(T-N)}\nu_{N}\left( \sum_{n=0}^{\infty} \frac{\lambda^n T^{n}}{n!} P_{final}^{n} \left( \int_{\frac{N}{T}}^{1} \phi(x)dx \right)^{n} \right) + \mathcal{E},
\end{align*}
where $\mathcal{E}$ denotes the rest of the terms. 
\vskip 0.2 in
\noindent
Take $K>0$ and define 
$$ T = \left(  \int_{K-1 \over K}^{1} \phi(x)dx \right)^{-1} t_{mix}(\epsilon/2) $$
and $$ N = {(K-1) \over K} \left(  \int_{K-1 \over K}^{1} \phi(x)dx \right)^{-1} t_{mix}(\epsilon/2)$$
\vskip 0.2 in
\noindent
Recall the approximation of the minimum function $\phi$ around $x=1$
$$\phi(x) = 1 +\frac{\phi^{(m)}(1)(x-1)^{m}}{m!}+O\left(|x-1|^{m+1}\right)$$
and therefore
$$ \int_{K-1 \over K}^{1} \phi(x)dx ={1 \over K}\left(1+{\gamma(K) \over K^m} \right),$$
where $\gamma(K)=(-1)^m{\phi^{(m)}(1) \over (m+1)!}+O(K^{-1})$.
Thus we can write
$$ \nu P_{T}(0,T) = e^{-\lambda(T-N)}\nu_{N}\left( \sum_{n=0}^{\infty} \frac{\lambda^n (T-N)^{n}}{n!} P_{final}^{n} \left[ 1 +\gamma(K) \left( \frac{T-N}{T} \right)^{m} \right]^{n} \right) + \mathcal{E}$$ 

\noindent
We see, using standard uniformization argument, that 
\begin{align*}
 \nu P_{T}(0,T) 
&= e^{-\lambda \left( 1+ {\gamma(K) \over K^m} \right)^{-1} t_{mix}(\epsilon/2)}\nu_{N}\left( \sum_{n=0}^{\infty} \frac{\lambda^n t_{mix}(\epsilon/2)^{n}}{n!} P_{final}^{n} \right) + \mathcal{E} \\
&= e^{\lambda \left({\gamma(K) \over K^m+\gamma(K)} \right) t_{mix}(\epsilon/2)}\nu_{N} \exp\left\{Q_{final}\cdot  t_{mix}(\epsilon/2)\right\}+ \mathcal{E}
\end{align*}

\noindent 
Now, since $(-1)^m\phi^{(m)}(1) \leq 0$, we have that, for any probability distribution $\nu$, 
$$\| \nu P_{T}(0,T) - \pi_{f} \|_{TV} =  \| \nu \exp\left\{Q_{final}\cdot t_{mix}(\epsilon/2)\right\} - \pi_{f} \|_{TV} \cdot e^{\lambda \left({\gamma(K) \over K^m+\gamma(K)} \right) t_{mix}(\epsilon/2)} + S_{N},$$
where, by triangle inequality,
$$0 \leq S_{N} \leq 1 - e^{\lambda \left( 1+ {\gamma(K) \over K^m} \right)^{-1} \cdot t_{mix}(\epsilon/2)} \left( \sum_{n=0}^{\infty} \frac{\lambda^{n} (t_{mix}(\epsilon/2))^{n}}{n!} \right)$$
and, by definition of $t_{mix}$,
$$ \| \nu \exp\left\{Q_{final}\cdot t_{mix}(\epsilon/2)\right\} - \pi_{f} \|_{TV} \cdot e^{\lambda \left({\gamma(K) \over K^m+\gamma(K)} \right) t_{mix}(\epsilon/2)} < \epsilon/2$$ 

\noindent
Taking $K=c\left(\lambda/\epsilon\right)^{1 \over m} t_{mix}(\epsilon/2)^{1 \over m}$ with constant $c>>(-1)^{m+1}{\phi^{(m)}(1) \over (m+1)!}$, we obtain
$$\epsilon>-\log(1 - \epsilon/2) \geq \lambda \left({-\gamma(K) \over K^m+\gamma(K)} \right) t_{mix}(\epsilon/2)$$ 
and therefore
$$0 \leq S_{N} \leq 1 - e^{\lambda \left({\gamma(K) \over K^m+\gamma(K)} \right) t_{mix}(\epsilon/2)} < \epsilon/2$$
\noindent
Thus
$$T={K t_{mix}(\epsilon/2) \over 1+{\gamma(K) \over K^m}}= \mathcal{O} \left( \left[{\lambda \over \epsilon}\right]^{1 \over m} t_{mix}^{\frac{m+1}{m}}(\epsilon/2) \right)$$
\end{proof}

\bibliographystyle{amsplain}

\end{document}